\documentclass[11pt]{amsart}
\usepackage{amsfonts} 
\usepackage{amsmath} 
\usepackage{amsthm} 
\usepackage{amssymb} 

\newcommand{\sep}{\itemsep 0pt}
\newcommand{\set}[2]{\{#1 \ ; \ #2\}}

\newcommand{\gen}[1]{\langle #1 \rangle}

\newcommand{\NN}{\mathbb{N}}
\newcommand{\ZZ}{\mathbb{Z}}
\newcommand{\QQ}{\mathbb{Q}}
\newcommand{\CC}{\mathbb{C}}

\newtheorem{theorem}{Theorem}[section]
\newtheorem{proposition}[theorem]{Proposition}
\newtheorem{lemma}[theorem]{Lemma}
\newtheorem{corollary}[theorem]{Corollary}

\theoremstyle{definition}
\newtheorem{definition}[theorem]{Definition}
\newtheorem{example}[theorem]{Example}
\newtheorem{remark}[theorem]{Remark}

\begin{document}

\title[Nilpotent, algebraic and quasi-regular elements]{Nilpotent, algebraic and quasi-regular elements in rings and algebras}
\author{Nik Stopar}
\address{Faculty of Electrical Engineering, University of Ljubljana, Tr\v za\v ska cesta 25, 1000 Ljubljana, Slovenia}
\email{nik.stopar@fe.uni-lj.si}

\begin{abstract}
We prove that an integral Jacobson radical ring is always nil, which extends a well known result from algebras over fields to rings. As a consequence we show that if every element $x$ of a ring $R$ is a zero of some polynomial $p_x$ with integer coefficients, such that $p_x(1)=1$, then $R$ is a nil ring. With these results we are able to give new characterizations of the upper nilradical of a ring and a new class of rings that satisfy the K\"{o}the conjecture, namely the integral rings.
\end{abstract}

\maketitle



\section{Introduction}

Let $R$ be an associative ring or algebra. Every nilpotent element of $R$ is quasi-regular and algebraic. In addition the quasi-inverse of a nilpotent element is a polynomial in this element. In the first part of this paper we will be interested in the connections between these three notions; nilpotency, algebraicity, and quasi-regularity. In particular we will investigate how close are algebraic elements to being nilpotent and how close are quasi-regular elements to being nilpotent. We are motivated by the following two questions:

\medskip\noindent\textbf{Q1.} Algebraic rings and algebras are usually thought of as nice and well behaved. For example an algebraic algebra over a field, which has no zero divisors, is a division algebra. On the other hand nil rings and algebras, which are of course algebraic, are bad and hard to deal with. It is thus natural to ask what makes the nil rings and algebras bad among all the algebraic ones.

\medskip The answer for algebras over fields is well known, namely they are Jacobson radical. We generalize this to rings (and more generally to algebras over Jacobson rings) in two different ways; firstly we show that nil rings are precisely those that are integral and Jacobson radical (see Theorem~\ref{ringI}), and secondly, we show that the only condition needed for an algebraic ring to be nil, is that its elements are zeros of polynomials $p$ with $p(1)=1$ (see Theorem~\ref{ring}). 

\medskip\noindent\textbf{Q2.} Can nilpotent elements be characterized by the property "quasi-inverse of $a$ is a polynomial in $a$"?

\medskip It is somewhat obvious that element-by-element this will not be possible, however we are able to characterize the upper nilradical in this way (see Corollary~\ref{nilradical}).

One of the most important problems concerning nil rings is the K\"{o}the conjecture. In 1930 K\"{o}the conjectured that if a ring has no nonzero nil ideals, then it has no nonzero nil one-sided ideals. The question whether this is true is still open. There are many statements that are equivalent to the K\"{o}the conjecture and many classes of rings and algebras that are known to satisfy the K\"{o}the conjecture (see \cite{Kre}, \cite{Smo2}, and \cite{Smo3} for an overview). We give yet another class of such rings, namely the integral rings (see Corollary~\ref{kothe}).

In the second part of this paper we investigate the structure of certain sets of elements of rings and algebras. In particular we show that a subgroup of the group of quasi-regular elements (equipped with quasi-multiplication) is closed for ring addition if and only if it is closed for ring multiplication. This gives us some information on the structure of the set of all elements of a ring which are zeros of polynomials $p$ with $p(1)=1$.

\section{Preliminaries}
Throughout this paper we are dealing with associative rings and algebras, possibly nonunital and noncommutative. Given a ring or algebra $(R,+,\cdot)$, we define an operation $\circ$ on $R$, called \emph{quasi-multiplication}, by
$$a \circ b=a+b-ab.$$
It is easy to see that $(R,\circ)$ is a monoid with identity element $0$. An element $a \in R$ is called \textit{quasi-regular} if it is invertible in $(R,\circ)$, i.e. if there exists $a' \in R$ such that $a \circ a'=a' \circ a=0$. In this case we say that $a'$ is the \textit{quasi-inverse} of $a$. If $R$ is unital then this is equivalent to $1-a$ being invertible in $(R,\cdot)$ with inverse $1-a'$. In fact the map $f: (R,\circ) \to (R,\cdot)$ given by $x \mapsto 1-x$ is a monoid homomorphism, since $1-a \circ b=(1-a)(1-b)$. The set of all quasi-regular elements of $R$ will be denoted by $Q(R)$. Clearly $(Q(R),\circ)$ is a group, since this is just the group of invertible elements of the monoid $(R,\circ)$. For every $a \in Q(R)$ and every $n \in \ZZ$ the $n$-th power of $a$ in $(Q(R),\circ)$ will be denoted by $a^{(n)}$ to distinguish it from $a^n$, the $n$-th power of $a$ in $(R,\cdot)$. In particular $a^{(0)}=0$ and $a^{(-1)}$ is the quasi-inverse of $a$. If $R$ is unital, then $1-a^{(-1)}=(1-a)^{-1}$. A subset $S \subseteq R$ is called \textit{quasi-regular} if $S \subseteq Q(R)$. The \textit{Jacobson radical} of $R$ is the largest quasi-regular ideal of $R$ and will be denoted by $J(R)$.

The set of all nilpotent elements in $R$ will be denoted by $N(R)$. Every nilpotent element is quasi-regular, so $N(R) \subseteq Q(R)$. In fact if $x^n=0$ then $-x-x^2-\ldots-x^{n-1}$ is the quasi-inverse of $x$. A subset $S \subseteq R$ is called \textit{nil} if $S \subseteq N(R)$. The \textit{upper nilradical} of $R$ is the largest nil ideal of $R$ and will be denoted by $Nil^\ast(R)$. If $R$ is commutative then $Nil^\ast(R)=N(R)$.

The \textit{lower nilradical} of $R$ (also called the \textit{prime radical}) is the intersection of all prime ideals of $R$ and will be denoted by $Nil_\ast(R)$. It can also be characterized as the lower radical determined by the class of all nilpotent rings (see \cite{Gar} for details). For any ring $R$ we have $Nil_\ast(R) \subseteq Nil^\ast(R) \subseteq J(R)$.

Let $K$ be a commutative unital ring and $R$ a $K$-algebra, possibly noncommutative and nonunital. An element $a \in R$ is \emph{algebraic} over $K$ if there exists a nonzero polynomial $p \in K[x]$ such that $p(0)=0$ and $p(a)=0$. If in addition $p$ can be chosen monic (i.e. the leading coefficient of $p$ is equal to $1$), then $a$ is called \emph{integral} over $K$. The condition $p(0)=0$ is necessary only because $R$ may be nonunital, in which case only polynomials with zero constant term can be evaluated at elements of $R$. The set of all algebraic elements of $R$ will be denoted by $A_K(R)$, the set of all integral elements of $R$ will be denoted by $I_K(R)$. A $K$-algebra $R$ is \emph{algebraic} (resp. \emph{integral}) over $K$ if every element in $R$ is algebraic (resp. integral) over $K$. Note the special case of the above definitions when $R$ is just a ring, in which case we considder it as an algebra over $K=\ZZ$. In this case we will also write $A(R)=A_\ZZ(R)$ and $I(R)=I_\ZZ(R)$. Clearly, every nilpotent element of $R$ is integral, so $N(R) \subseteq I_K(R) \subseteq A_K(R)$. If $K=F$ is a field then $I_F(R)=A_F(R)$.

\section{$\pi$-algebraic rings and algebras}

Throughout this section $K$ will always denote a commutative unital ring, $F$ a field, and $R$ an algebra over $K$ or $F$, unless specified otherwise. The two questions from the introduction motivate the following definition, which will play a crucial role in our considerations.

\begin{definition}
An element $a$ of a $K$-algebra $R$ is \textit{$\pi$-algebraic} (over $K$) if there exists a polynomial $p \in K[x]$ such that $p(0)=0$, $p(1)=1$ and $p(a)=0$. In this case we will also say that $a$ is \textit{$\pi$-algebraic with polynomial $p$}. A subset $S \subseteq R$ is \textit{$\pi$-algebraic} if every element in $S$ is $\pi$-algebraic. The set of all $\pi$-algebraic elements of a $K$-algebra $R$ will be denoted by $\pi_K(R)$.
\end{definition}

When $R$ is just a ring, we considder it as an algebra over $K=\ZZ$, and write $\pi(R)=\pi_\ZZ(R)$. The crucial condition in this definition is the condition $p(1)=1$. The condition $p(0)=0$ is there simply because $R$ may not be unital, in which case only polynomials with zero constant term can be evaluated at an element of $R$.

We first present some basic properties of $\pi$-algebraic elements along with some examples.

\begin{lemma}\label{inclusion}
If $R$ is a $K$-algebra then $N(R) \subseteq \pi_K(R) \subseteq A_K(R) \cap Q(R)$. If $R$ is an $F$-algebra then $N(R) \subseteq \pi_F(R)=A_F(R) \cap Q(R)$. The quasi-inverse of a $\pi$-algebraic element is a polynomial in this element.
\end{lemma}

\proof Clearly every nilpotent element is $\pi$-algebraic and every $\pi$-algebraic element is algebraic. Suppose $a \in R$ is $\pi$-algebraic with polynomial $p$. Then $P(x)=1-(1-p(x))/(1-x)$ is again a polynomial with $P(0)=0$ (and proper coefficients). Hence we may define $a'=P(a)$. Since $x \circ P(x)=x+P(x)-xP(x)=p(x)$, we have $a \circ a'=0$. Similarly we get $a' \circ a=0$. Hence $a'$ in the quasi-inverse of $a$ and it is a polynomial in $a$. Now suppose $R$ is an $F$-algebra and $a$ is an element of $A_F(R) \cap Q(R)$. Let $r \in F[x]$ be the minimal polynomial of $a$ (if $R$ is not unital then $r(0)$ must be zero) and let $a'$ be the quasi-inverse of $a$. Suppose $r(1)=0$. Then $r(x)=(1-x)q(x)=q(x)-xq(x)$ for some polynomial $q \in F[x]$ of degree less then that of $r$. If $R$ is not unital then $q(0)=0$, so we may evaluate $q$ at $a$ in any case. Hence $0=r(a)-a'r(a)=q(a)-aq(a)-a'q(a)+a'aq(a)=q(a)-(a' \circ a)q(a)=q(a)$, which is a contradiction since $r$ was the minimal polynomial for $a$. Thus $r(1)$ is an invertible element of $F$ and hence the element $a$ is $\pi$-algebraic with polynomial $p(x)=r(1)^{-1}r(x)x$. \endproof

We shall see in the examples that the inclusion $\pi_K(R) \subseteq A_K(R) \cap Q(R)$ may be strict.

\begin{lemma}\label{2minus}
If $R$ is a unital $K$-algebra then $2-\pi_K(R) \subseteq \pi_K(R)$. In particular $0,2 \in \pi_K(R)$ and $1 \notin \pi_K(R)$. If $R$ is a unital $F$-algebra then $2-\pi_F(R) \subseteq \pi_F(R)$. In addition $F\backslash \{1\} \subseteq \pi_F(R)$ and $1 \notin \pi_F(R)$.
\end{lemma}

\proof If $a$ is $\pi$-algebraic with polynomial $p$ then $2-a$ is $\pi$-algebraic with polynomial $q(x)=p(2-x)x$. We always have $0 \in \pi_K(R)$, hence $2 \in \pi_K(R)$. The identity element is never $\pi$-algebraic since it is not quasi-regular. If $R$ is a unital $F$-algebra and $\lambda \neq 1$ is a scalar then $\lambda$ is $\pi$-algebraic with polynomial $p(x)=(1-\lambda)^{-1}(x-\lambda)x$. \endproof

Next we give a few examples.

\begin{example}
For a finite ring $R$, $\pi(R)=Q(R)$ and $J(R)=Nil^\ast(R)$. To verify the first part observe that $(Q(R),\circ)$ is a finite group, say of order $n$. So for every $a \in Q(R)$ we have $a^{(n)}=0$, hence every $a \in Q(R)$ is $\pi$-algebraic with polynomial $p(x)=x^{(n)}=1-(1-x)^n$. The second part is well known and it also follows from the first part and Theorem~\ref{ring}.
\end{example}

\begin{example}\label{rac}
For any field $F$, $\pi_F(F)=F\backslash \{1\}=Q(F)$ by Lemma~\ref{2minus}. In particular $\pi_{\QQ}(\QQ)=\QQ\backslash \{1\}=Q(\QQ)$. On the other hand we have $\pi(\QQ)=\set{1+\frac{1}{n}}{n \in \ZZ\backslash \{0\}}$. Indeed, if $n$ is a nonzero integer then $1+\frac{1}{n}$ is $\pi$-algebraic over $\ZZ$ with polynomial $s(x)=(1-n(x-1))x$. Conversely, suppose $\frac{a}{b} \in \QQ$ with $a$ and $b$ coprime, is $\pi$-algebraic with polynomial $p \in \ZZ[x]$ of degree $d$. Then $q(x)=b^dp(\frac{x}{b})$ is a polynomial with integer coefficients. Hence $a-b$ divides $q(a)-q(b)=b^dp(\frac{a}{b})-b^dp(1)=-b^d$. Since $a$ and $b$ are coprime, this is only possible if $a-b=\pm 1$ (any prime that would divide $a-b$ would divide $b$ and hence $a$). Thus $\frac{a}{b}=1\pm \frac{1}{b}$ as needed. Obviously $A(\QQ)=\QQ$, so the inclusion $\pi(\QQ) \subseteq A(\QQ) \cap Q(\QQ)$ from Lemma~\ref{inclusion} is strict here.
\end{example}

\begin{example}\label{mat}
Let $F \subseteq E$ be fields and $M_n(E)$ the ring of $n \times n$ matrices over $E$. Then
\begin{eqnarray*}
N(M_n(E)) &=& \textrm{ matrices with eigenvalues $0$}, \\
\pi_F(M_n(E)) &=& \textrm{ matrices with eigenvalues in $\overline{F}\backslash \{1\}$}, \\
Q(M_n(E)) &=& \textrm{ matrices with eigenvalues in $\overline{E}\backslash \{1\}$},
\end{eqnarray*}
where $\overline{F}\subseteq\overline{E}$ are algebraic closures of $F$ and $E$.
A matrix is quasi-regular iff it has no eigenvalue equal to $1$. So in view of Lemma~\ref{inclusion}, to verify the above, we only need to prove that
$$A_F(M_n(E))= \textrm{ matrices with eigenvalues in $\overline{F}$}.$$
If $A \in M_n(E)$ is algebraic over $F$, it clearly has eigenvalues in $\overline{F}$. So suppose $A \in M_n(E)$ has eigenvalues $\lambda_1,\lambda_2,\ldots,\lambda_n \in \overline{F}$. For every $i=1,2,\ldots,n$, let $p_i$ be the minimal polynomial of $\lambda_i$ over $F$. Then the minimal polynomial $m_A$ of $A$ over $E$ divides $P(x)=\prod_{i=1}^{n}p_i(x)$, hence $P(A)=0$. Since $P$ has coefficients in $F$, $A$ is algebraic over $F$.
\end{example}

The following proposition gives a connection between $\pi$-algebraic and integral elements.

\begin{proposition}\label{integral}
Let $R$ be a $K$-algebra and $a$ an element of $R$. The following are equivalent:
\begin{enumerate}\sep
\item\label{int1} $a$ is $\pi$-algebraic,
\item\label{int2} $a$ is quasi-regular and $a^{(-1)}$ is integral,
\item\label{int3} $a$ is quasi-regular and $a^{(-1)}$ is a polynomial in $a$.
\end{enumerate}
\end{proposition}

\proof By Lemma~\ref{inclusion} \eqref{int1} implies \eqref{int3}. On the other hand, if $a^{(-1)}=P(a)$ where $P$ is a polynomial in $K[x]$, then $a+P(a)-aP(a)=0$, so $a$ is $\pi$-algebraic with polynomial $(x+P(x)-xP(x))x$. It remains to prove the equivalence of \eqref{int1} and \eqref{int2}.

For a polynomial $p \in K[x]$ define $\widehat{p}(x)=(x-1)^{\deg p}p(\frac{x}{x-1})$, which is again a polynomial in $K[x]$. Notice that $\widehat{p}(1)$ equals the leading coefficient of $p$ and the leading coefficient of $\widehat{p}$ equals $p(1)$ (the sum of all coefficients of $p$) if $p(1)\neq 0$. In addition $\widehat{p}(0)=0$ iff $p(0)=0$.

We may assume that $R$ is unital, otherwise we just adjoin a unit to $R$. Let $a$ be a quasi-regular element. Then the inverse of $1-a$ is $1-a^{(-1)}$, so the term $\frac{x}{x-1}$ evaluated at $a$ equals $-a(1-a^{(-1)})=-a+aa^{(-1)}=a^{(-1)}$. Thus $\widehat{p}(a)=(a-1)^{\deg p}p(a^{(-1)})$. This shows that $\widehat{p}(a)=0$ iff $p(a^{(-1)})=0$, since $1-a$ is invertible. Similarly $\widehat{p}(a^{(-1)})=0$ iff $p(a)=0$.

If $p$ is a monic polynomial such that $p(0)=0$ and $p(a^{(-1)})=0$ then $a$ is $\pi$-algebraic with polynomial $\widehat{p}$. If $a$ is $\pi$-algebraic with polynomial $p$ then $\widehat{p}$ is a monic polynomial such that $\widehat{p}(0)=0$ and $\widehat{p}(a^{(-1)})=0$, so $a^{(-1)}$ is integral. \endproof

In particular Proposition~\ref{integral} states that $\pi_K(R)=(Q(R) \cap I_K(R))^{(-1)}$ (compare with Lemma~\ref{inclusion}).

By Lemma~\ref{inclusion} an algebra over a field $F$ is $\pi$-algebraic if and only if it is algebraic and Jacobson radical. So the following proposition is just a restatement of a well known fact that any algebraic Jacobson radical $F$-algebra is nil (see for example \cite[p. 144]{Sza}). In fact, every algebraic element in the Jacobson radical of an $F$-algebra is nilpotent and its nilindex is equal to its degree.

\begin{proposition}\label{alg}
Every $\pi$-algebraic $F$-algebra is nil.
\end{proposition}

We now extend this result to algebras over Jacobson rings. Recall that a commutative unital ring $K$ is a \emph{Jacobson ring} (or a \emph{Hilbert ring}) if every prime ideal of $K$ is an intersection of maximal ideals of $K$. Examples of Jacobson rings are fields and polynomial rings over fields in finitely many commutative variables. In addition, any principal ideal domain with infinitely many irreducible elements is also a Jacobson ring. In particular, the ring of integers $\ZZ$ is a Jacobson ring.

\begin{theorem}\label{ringI}
If $K$ is a Jacobson ring then every integral Jacobson radical $K$-algebra is nil.
\end{theorem}

\proof Let $R$ be an integral Jacobson radical $K$ algebra and $a \in R$. Consider $R$ as a subalgebra of some unital $K$-algebra $R^1$. Let $K[a]$ be a unital subalgebra of $R^1$ generated by $a$. Since $K$ is a Jacobson ring and $K[a]$ is a finitely generated (commutative unital) $K$-algebra, $K[a]$ is a Jacobson ring by a version of Hilbert's Nullstellensatz \cite[Theorem~4.19]{Eis}. Hence $J(K[a])=Nil_\ast(K[a])$ is a nil ideal. It suffices to prove that $a \in J(K[a])$. Take any $r \in K[a]$. Since $ar$ is an element of $R$, it is quasi-regular in $R$ and its quasi-inverse $(ar)^{(-1)} \in R$ is integral. By Proposition~\ref{integral}, $(ar)^{(-1)}$ is a polynomial in $ar$. But $ar$ is a polynomial in $a$, hence $(ar)^{(-1)} \in K[a]$, i.e. $ar$ is quasi-regular in $K[a]$. Since $r$ was arbitrary, we conclude that $a \in J(K[a])$. \endproof

Without the assumption that the ring $K$ is Jacobson, Theorem~\ref{ringI} fails.

\begin{proposition}\label{ringI-converse}
If $K$ is not a Jacobson ring then there exists an integral Jacobson radical $K$-algebra which is not nil.
\end{proposition}

\proof Let $P$ be a prime ideal of $K$ which is not an intersection of maximal ideals. Then $J(K/P)$ is a nonzero $K$-algebra. In addition, it is Jacobson radical and integral because an element $k+P \in J(K/P)$ is integral over $K$ with polynomial $x^2-kx$. Since $K$ is commutative and $P$ is a prime ideal, the algebra $K/P$ has no nonzero nilpotent elements, hence $J(K/P)$ is not nil. \endproof

The assumption that the algebra is integral in Theorem~\ref{ringI} is also crucial. A merely algebraic Jacobson radical algebra over a Jacobson ring need not be nil.

\begin{example}
Consider $R=\set{\frac{2m}{2n-1}}{m,n \in \ZZ}$ as a subring of rational numbers. The quasi-inverse of $\frac{2m}{2n-1}$ is $\frac{2m}{2m-2n+1}$, which is again an element of $R$. So $R$ is a Jacobson radical ring algebraic over $\ZZ$, but it is not nil.
\end{example}

As a direct consequence of Theorem~\ref{ringI} and Proposition~\ref{integral} we get

\begin{theorem}\label{ring}
If $K$ is a Jacobson ring then every $\pi$-algebraic $K$-algebra is nil.
\end{theorem}

This answers question Q1 in two ways: the fact that distinguishes nil rings and algebras from all other algebraic ones is firstly that they are integral and Jacobson radical, and secondly that the polynomials ensuring algebraicity in the nil case have the sum of their coefficients equal to $1$. It is perhaps interesting that this rather large family of polynomials with the sum of coefficients equal to $1$ produces the same effect as the rather restrictive family $\{x,x^2,x^3,x^4,\ldots\}$.

Observe that in an algebraic division $F$-algebra only the identity is not $\pi$-algebraic (since all other elements are quasi-regular). So if only one element in an algebra is not $\pi$-algebraic then the algebra may be very nice instead of nil.

Next corollary addresses question Q2, giving new characterizations of the upper nilradical of a ring in the process. We formulate it only for rings, though it is valid for all algebras over Jacobson rings.

\begin{corollary}\label{nilradical}
For a ring $R$ the following hold:
\begin{enumerate}\sep
\item $Nil^\ast(R)$ is the largest $\pi$-algebraic ideal of $R$,
\item $Nil^\ast(R)$ is the largest integral quasi-regular ideal of $R$,
\item $Nil^\ast(R)$ is the largest quasi-regular ideal of $R$ such that the quasi-inverse of each element is a polynomial in this element.
\end{enumerate}
\end{corollary}

\proof If $I$ is an ideal of $R$ satisfying any of the above conditions then $I$ is $\pi$-algebraic by Proposition~\ref{integral} and thus nil by Theorem~\ref{ring}. Hence $Nil^\ast(R)$ is the largest such ideal. \endproof

\begin{corollary}\label{JN}
If $R$ is an integral ring then $J(R)=Nil^\ast(R)$.
\end{corollary}

A ring $R$ is said to satisfy the K\"{o}the conjecture if every nil one-sided ideal of $R$ is contained in a nil two-sided ideal of $R$ (cf. \cite{Yon}). If $J(R)=Nil^\ast(R)$ for a ring $R$, then $R$ satisfies the K\"{o}the conjecture since $J(R)$ contains every nil one-sided ideal. Corollary~\ref{JN} thus implies

\begin{corollary}\label{kothe}
Every integral ring satisfies the K\"{o}the conjecture.
\end{corollary}

In what follows we will exhibit an even stronger connection between $\pi$-algebraic and nilpotent elements than that given by Theorem~\ref{ring}, in case the ring $K$ satisfies certain properties given by the following definition.

\begin{definition}
We shall say that a principal ideal domain $K$ is \emph{exceptional} if there is no nonconstant polynomial $p \in K[x]$ such that $p(k)$ would be invertible in $K$ for all $k \in K$.
\end{definition}

Exceptional PIDs are quite common, here are some examples.

\begin{proposition}\label{exceptional}
\mbox{}
\begin{enumerate}\sep
\item\label{ex1} A field is an exceptional PID if and only if it is algebraically closed.
\item\label{ex2} The ring of integers $\ZZ$ and the ring of Gaussian integers $\ZZ[i]$ are exceptional PIDs.
\item\label{ex3} For any field $F$ the polynomial ring $F[x]$ is an exceptional PID.
\item\label{ex4} If $K$ is an exceptional PID and $S \subseteq K$ is a multiplicatively closed subset multiplicatively generated by a finite number of elements, then the localization $S^{-1}K$ is an exceptional PID.
\end{enumerate}
\end{proposition}

\proof Claim \eqref{ex1} is clear.

\smallskip\noindent \eqref{ex2}: Let $p$ be a polynomial in $\ZZ[x]$ such that $p(k)$ is invertible for all $k \in \ZZ$. Since there are only finitely many invertible elements in $\ZZ$, there exist an invertible element $u \in \ZZ$, such that $p(k)=u$ for infinitely many $k \in \ZZ$. But then the polynomial $p(x)-u$ has infinitely many zeros, so it must be zero. Hence $p$ is a constant polynomial. The same proof works for $\ZZ[i]$.

\smallskip\noindent \eqref{ex3}: Let $F$ be a field and $P(y)$ a nonconstant polynomial in $(F[x])[y]$. Write $P(y)=p_0(x)+p_1(x)y+\ldots+p_n(x)y^n$ where $p_n(x) \neq 0$ and $n \geq 1$. Denote $d_i=\deg p_i$ and $d=\max\set{d_i}{i=0,1,2,\ldots,n}$, where the degree of the zero polynomial is equal to $-\infty$. Let $p(x)=x^{d+1}$. The degree of $p_i(x)(p(x))^i$ is equal to $d_i+i(d+1)$. Since $d_n,d \neq -\infty$, we have $d_n+n(d+1) \geq n(d+1)>d+(n-1)(d+1) \geq d_i+i(d+1)$ for all $i<n$. This implies that the degree of $P(p(x))$ is equal to $d_n+n(d+1) \geq 1$, hence $P(p(x))$ is not invertible in $F[x]$.

\smallskip\noindent \eqref{ex4}: A localization of a PID is again a PID. Factor each generator of $S$ into irreducible factors and let $S' \subseteq K$ be a multiplicatively closed subset multiplicatively generated by all the irreducible elements appearing in these factorizations. Since the localizations of $K$ at $S$ and $S'$ are isomorphic, we may assume that $S=S'$, i.e. $S$ is generated by a finite number of irreducible elements. Now suppose $p(x) \in S^{-1}K[x]$ is a polynomial, such that $p(\widehat{k})$ is invertible in $S^{-1}K$ for all $\widehat{k} \in S^{-1}K$. Take any $s \in S$, such that the coefficients of $sp(x)$ are elements of $K$. Let $t$ be the product of all irreducible elements in $S$. Observe that the coefficients of the polynomial $sp(sp(0)tx)$ are elements of $K$ and are all divisible by $sp(0) \in K$. Hence, $P(x)=\frac{s}{sp(0)}p(sp(0)tx)=\frac{1}{p(0)}p(sp(0)tx)$ is a polynomial with coefficients in $K$. Now take any $k \in K$. By assumption, $p(0)$ and $p(sp(0)tk)=p(0)P(k)$ are invertible in $S^{-1}K$, hence so is $P(k)$. But $P(k) \in K$, so the only irreducible elements that may divide $P(k)$ are those that lie in $S$. However, any such irreducible element divides $t$, hence it divides all coefficients of $P$ except $P(0)=1$, so it cannot divide $P(k)$. This shows that $P(k)$ is invertible in $K$. Since $K$ is an exceptional PID, $P(x)$, and consequently $p(x)$, must be a constant polynomial. \endproof

In a PID every nonzero prime ideal is maximal, so a PID is a Jacobson ring if and only if $0$ is an intersection of maximal ideals, i.e. the Jacobson radical is $0$.

\begin{proposition}\label{J-semisimple}
If $K$ is an exceptional PID then $J(K)=0$, i.e. $K$ is a Jacobson ring. In particular, if $K$ is not a field then $K$ has infinitely many nonassociated irreducible elements.
\end{proposition}

\proof Let $K$ be an exceptional PID. Suppose $J(K) \neq 0$ and take $0 \neq a \in J(K)$. Since $K$ is commutative and unital, this implies that $1-ak$ is invertible in $K$ for every $k \in K$. But then the polynomial $p(x)=1-ax$ contradicts the definition of an exceptional PID. Hence $J(K)=0$. Since $K$ is commutative and unital, $J(K)$ is just the intersection of all maximal ideals of $K$. If $K$ is not a field then the maximal ideals of $K$ are the principal ideals generated by the irreducible elements. If there are only finitely many such ideals then their intersection is nonzero. \endproof

The converse of Proposition~\ref{J-semisimple} does not hold. There exist PIDs which are Jacobson rings but are not exceptional. The simplest example is given by any field that is not algebraically closed, however fields are rather extremal among all PID, since they have no irreducible elements. Hence, we give an example which is not a field.

\begin{example}
Let $S \subseteq \ZZ$ be a multiplicatively closed subset multiplicatively generated by all primes $p$ with $p=2$ or $p \equiv 1 \pmod 4$ and let $K=S^{-1}\ZZ$ be the localization of $\ZZ$ at $S$. Then $K$ has infinitely many nonassociated irreducible elements, represented by the primes $p$ with $p \equiv 3 \pmod 4$, hence $J(K)=0$ and $K$ is a Jacobson ring. Now let $p(x)=x^2+1$. To see that $K$ is not exceptional, we will show that $p(k)$ is invertible in $K$ for all $k \in K$. For $k=\frac{m}{n} \in K$ we have $p(k)=\frac{m^2+n^2}{n^2}$. To see that this is invertible in $K$ we need to show that any prime dividing $m^2+n^2$ is contained in $S$. Suppose $p$ is a prime with $p \equiv 3 \pmod 4$ that divides $m^2+n^2$. Then $m^2 \equiv -n^2 \pmod p$. Since $n \in S$, this implies that both $m$ and $n$ are coprime to $p$. Hence we have $1 \equiv m^{p-1} \equiv (m^2)^{\frac{p-1}{2}} \equiv (-n^2)^{\frac{p-1}{2}} \equiv (-1)^{\frac{p-1}{2}}n^{p-1} \equiv (-1)^{\frac{p-1}{2}} \equiv -1 \pmod p$. This is a contradiction since $p \neq 2$, which finishes the proof.
\end{example}

Theorem~\ref{ring} implies that if the subalgebra generated by an element $a$ is $\pi$-algebrac, then $a$ is a nilpotent element. The next proposition, which was our main motivation for the introduction of exceptional PIDs, considers the situation when only the submodule generated by $a$ is assumed to be $\pi$-algebraic. It thus gives a stronger connection between $\pi$-algebraic and nilpotent elements for algebras over exceptional PIDs.

\begin{proposition}\label{na}
Let $K$ be an exceptional principal ideal domain and $R$ a $K$-algebra. If $a$ is an element of $R$ such that $Ka \subseteq \pi_K(R)$, then there exists $0 \neq k \in K$ such that $ka$ is nilpotent. In particular, if $R$ has no $K$-torsion, then $a$ is nilpotent.
\end{proposition}

\proof For a nonzero polynomial $f \in K[x]$, let $\delta(f)$ denote the greatest common divisor of all coefficients of $f$. First we show that for any $\pi$-algebraic element $r$ there exists a nonzero polynomial $f \in K[x]$ and a nonzero element $c \in K$ such that $f(1)=1$, $cf(r)=0$, and $f$ divides (within $K[x]$) any polynomial that annihilates $r$. So let $r$ be $\pi$-algebraic with polynomial $h \in K[x]$. Choose a nonzero polynomial $p \in K[x]$ of minimal degree such that $p(r)=0$ and let $c=\delta(p)$ and $f(x)=\frac{p(x)}{c} \in K[x]$. So $cf(r)=0$ and $\delta(f)=1$. Suppose $P \in K[x]$ is a polynomial with $P(r)=0$. By the division algorithm there exists $0 \neq m \in K$ and polynomials $s,t \in K[x]$ with $\deg t<\deg f=\deg p$ such that $mP(x)=s(x)f(x)+t(x)$ (divide in $F[x]$, where $F$ is the field of fractions of $K$, and multiply by a common denominator of all fractions). Multiplying by $c$ we get $cmP(x)=cs(x)f(x)+ct(x)=s(x)p(x)+ct(x)$. The minimality of $p$ now implies $ct(x)=0$, hence $t(x)=0$ and $mP(x)=s(x)f(x)$. By Gauss's lemma this implies $\delta(s)=m\delta(P)$ up to association, so $m$ divides $\delta(s)$. Thus the polynomial $\frac{s(x)}{m}$ has coefficients in $K$ and $P(x)=\frac{s(x)}{m}f(x)$, i.e. $f$ divides $P$. In particular $f$ divides $h$, so there is a polynomial $S$ such that $h(x)=S(x)f(x)$. Evaluating at $1$ we get $1=S(1)f(1)$, so $f(1)$ is invertible in $K$. We may assume that $f(1)=1$, otherwise we just multiply $f$ by $f(1)^{-1}=S(1)$.

Now let $R$ be a $K$-algebra and $a$ an element of $R$ with $Ka \subseteq \pi_K(R)$. By the above, for any $k \in K$ there exists $0 \neq c_k \in K$ and $0 \neq f_k \in K[x]$ such that $f_k(1)=1$, $c_kf_k(ka)=0$, and $f_k$ divides any polynomial that annihilates $ka$. Let $k \neq 0$. Then $f_1$ divides $c_kf_k(kx)$, since $c_kf_k(kx)$ annihilates $a$. Similarly $c_1k^{\deg f_1}f_1(\frac{x}{k})$ is a polynomial in $K[x]$ that annihilates $ka$, so $f_k$ divides $c_1k^{\deg f_1}f_1(\frac{x}{k})$. This in particular implies that all these polynomials have the same degree, so there exists $d_k \in K$ such that $c_1k^{\deg f_1}f_1(\frac{x}{k})=d_kf_k(x)$. We have $f_k(1)=1$, hence $\delta(f_k)=1$. Consequently, $c_1\delta(k^{\deg f_1}f_1(\frac{x}{k}))=d_k$ up to association. If $k$ is coprime to the leading coefficient of $f_1$ then $\delta(k^{\deg f_1}f_1(\frac{x}{k}))=1$ since $\delta(f_1)=1$. For such $k$ we have $c_1=d_k$ up to association, hence $c_1$ divides $d_k$ and $u_k=\frac{d_k}{c_1}$ is invertible. In addition $k^{\deg f_1}f_1(\frac{x}{k})=u_kf_k(x)$. Evaluating at $1$, we get $k^{\deg f_1}f_1(\frac{1}{k})=u_k$. Now $p(x)=x^{\deg f_1}f_1(\frac{1}{x})$ is a polynomial in $K[x]$ with $p(0)$ equal to the leading coefficient of $f_1$. Hence, we have proved above that $p(k)$ is invertible for every $k \neq 0$ coprime to $p(0)$. If we define $t(x)=p(p(0)x-1) \in K[x]$, then $t(k)$ is invertible for all $k \in K$ ($p(0)k-1=0$ means that $p(0)$ is invertible). Since $K$ is exceptional, it follows that $t$ is a constant polynomial and so is $p$. Hence, there exists $u \in K$ such that $f_1(\frac{1}{x})=\frac{u}{x^{\deg f_1}}$, i.e. $f_1(x)=ux^{\deg f_1}$. Consequently, $c_1ua^{\deg f_1}=0$ and $c_1ua$ is nilpotent. Clearly $c_1u \neq 0$. \endproof

\begin{remark}\label{rem-na}
We shall later need a slightly modified version of Proposition~\ref{na} with $K=\ZZ$. Observe that the conclusion still holds if we assume just $\NN a \subseteq \pi(R)$ instead of $\ZZ a \subseteq \pi(R)$. Indeed, one just has to replace polynomial $t(x)=p(p(0)x-1)$ in the proof with the polynomial $\widehat{t}(x)=p((p(0)x-1)^2)$.
\end{remark}

Without the assumption that $K$ is exceptional, Proposition~\ref{na} fails.

\begin{proposition}
Let $K$ be a principal ideal domain which is not exceptional. Then there exists a $K$-algebra $R$ and an element $a \in R$ such that $Ka \subseteq \pi_K(R)$, but $ka$ is not nilpotent for any $0 \neq k \in K$.
\end{proposition}

\proof Choose a nonconstant polynomial $p \in K[x]$, such that $p(k)$ is invertible in $K$ for all $k \in K$. Let $F$ be the algebraic closure of the field of fractions of $K$. Clearly $F$ is a $K$-algebra. Since polynomial $p$ is nonconstant, the polynomial $P(x)=x^{\deg p}p\left(\frac{1}{x}\right) \in K[x]$ has a nonzero root $a \in F$. Clearly, $ka$ is not nilpotent for any $0 \neq k \in K$. To finish the proof we show, that $Ka \subseteq \pi_K(F)$. The zero element is always $\pi$-algebraic, so take any $0 \neq k \in K$. Observe that $\deg P=\deg p$, because $p(0)$ is invertible. Hence $Q(x)=k^{\deg p}P\left(\frac{x}{k}\right) \in K[x]$. Since $Q(1)=p(k)$ is invertible in $K$, the element $ka$ is $\pi$-algebraic over $K$ with polynomial $Q(1)^{-1}Q(x)x \in K[x]$. \endproof

Recall that an algebra $R$ is called \textit{nil of bounded index $\leq n$} if $a^n=0$ for all $a \in R$. $R$ is called \textit{nil of bounded index} if there exists an integer $n$ such that $R$ is nil of bounded index $\leq n$. Similarly we will say that a $K$-algebra $R$ is \textit{$\pi$-algebraic of bounded degree $\leq n$} (resp. \textit{integral of bounded degree $\leq n$}) if every element of $R$ is $\pi$-algebraic (resp. integral) over $K$ with some polynomial of degree $\leq n$. $R$ is \textit{$\pi$-algebraic of bounded degree} (resp. \textit{integral of bounded degree}) if there exists an integer $n$ such that $R$ is $\pi$-algebraic of bounded degree $\leq n$ (resp. integral of bounded degree $\leq n$).

It follows from the proof of Proposition~\ref{integral} that an algebra is $\pi$-algebraic of bounded degree $\leq n$ if and only if it is Jacobson radical and integral of bounded degree $\leq n$. Theorem~\ref{ring} raises the following natural question. If an algebra $R$ over a Jacobson ring $K$ is $\pi$-algebraic of bounded degree, is it nil of bounded index? The answer is positive for algebras with no $K$-torsion.

\begin{corollary}\label{bounded}
Let $K$ be a Jacobson ring. If $R$ is a $\pi$-algebraic $K$-algebra of bound degree $\leq n$ with no $K$-torsion, then $R$ is nil of bounded index $\leq n$.
\end{corollary}

\proof Let $R$ be a $\pi$-algebraic $K$-algebra of bounded degree $\leq n$ with no $K$-torsion. By the remark above, $R$ is integral of bounded degree $\leq n$, and by Theorem~\ref{ring}, $R$ is nil. Take any $a \in R$. Let $p \in K[x]$ be a monic polynomial of degree $\leq n$, such that $p(a)=0$, and let $m$ be the smallest integer such that $a^m=0$. Suppose $m>n$. Write $p$ in the form $p(x)=t(x)x^k$, where $t(0) \neq 0$ and $k\leq n<m$. Multiplying the equality $0=t(a)a^k$ by $a^{m-k-1}$, we get $0=t(a)a^{m-1}=t(0)a^{m-1}$, because $a^m=0$. Since $R$ has no $K$-torsion, this implies $a^{m-1}=0$, which is in contradiction with the choice of $m$. Thus $m \leq n$ as needed. \endproof

Perhaps surprisingly, the answer for general algebras over Jacobson rings is negative as the following example shows.

\begin{example}
Let $K$ be a Jacobson PID, which is not a field. Then $K$ has infinitely many nonassociated irreducible elements. Choose a countable set of nonassociated irreducible elements $\{p_1,p_2,p_3,\ldots\}$ and let $R=\bigoplus_{i=1}^{\infty} p_iK/p_i^iK$. Clearly, $R$ is nil, but not of bounded index. Let $a=(a_i)_i$ be an element of $R$. By the Chinese remainder theorem there is an element $k \in K$ such that $k \equiv a_i \pmod{p_i^i}$ for all $i$ with $a_i \neq 0$. Thus $a$ is a zero of the monic polynomial $x^2-kx$. This shows that $R$ is integral of bounded degree $\leq 2$, hence it is also $\pi$-algebraic of bounded degree $\leq 2$.
\end{example}

Nevertheless the following holds for arbitrary algebras over Jacobson rings.

\begin{proposition}
Let $K$ be a Jacobson ring. If $R$ is a $\pi$-algebraic $K$-algebra of bounded degree then $Nil_\ast(R)=R$. In particular, $R$ is locally nilpotent.
\end{proposition}

\proof Suppose $P$ is a prime ideal of $R$. We want to apply Corollary~\ref{bounded} to $R/P$. $K$-algebra $R/P$ is again $\pi$-algebraic of bounded degree. Let $I=\set{k \in K}{k(R/P)=0}$. Clearly, $I$ is an ideal of $K$ and $R/P$ becomes a $K/I$-algebra if we define $(k+I)(r+P)=k(r+P)=kr+P$. Observe that $R/P$ is $\pi$-algebraic of bounded degree also over $K/I$. In addition, $R/P$ has no $K/I$-torsion. Indeed, if $(k+I)(r+P)=0$ for some $k \in K$ and $r \in R$ with $r+P \neq 0$, then $J=\set{x+P \in R/P}{k(x+P)=0}$ is a nonzero ideal of $R/P$. But $k(R/P) \cdot J=0$ and $R/P$ is a prime $K$-algebra, so $k(R/P)=0$, i.e. $k+I=0$ in $K/I$ as needed. $K/I$ is again a Jacobson ring, hence Corollary~\ref{bounded} implies that $R/P$ is nil of bounded index. Thus, by a result of Levitzki \cite[Theorem 4]{Lev}, we have $Nil_\ast(R/P)=R/P$, but on the other hand, $Nil_\ast(R/P)=0$ since $P$ is a prime ideal. So $P=R$, which shows that $Nil_\ast(R)=R$. \endproof

\section{The structure of $\pi(R)$}

In this section we investigate the structure of the set of all $\pi$-algebraic elements of an algebra. We restrict ourselves to algebras over fields and to rings. Throughout the section, $F$ will always denote a field and $R$ an $F$-algebra or a ring.

Recall that $(Q(R),\circ)$ is a group and by Lemma~\ref{inclusion} we have $N(R) \subseteq \pi(R) \subseteq Q(R)$. It is thus natural to ask under what conditions $N(R)$ and $\pi(R)$ are subgroups of $Q(R)$ and more generally what can be said about the structure of $\pi(R)$. In general $\pi(R)$ will not be closed under $\circ$. We give a concrete example later (see Example~\ref{notcl}), but the reason for this is that the integral elements of $R$ do not have any structure in general (they do not form a subring). However, if $R$ is commutative, then $\pi(R)$ will be closed under $\circ$. From here on $Q(R)$ will always be considered as a group with operation $\circ$. 

\begin{lemma}\label{aut}
For a quasi-regular element $r \in R$ the map $x \mapsto r\circ x\circ r^{(-1)}$ is an automorphism of $R$.
\end{lemma}

The proof of this lemma is an easy calculation. In fact, if $R$ is unital then $r\circ x\circ r^{(-1)}=(1-r)x(1-r^{(-1)})$, so the map is just the usual conjugation by $1-r$.

\begin{proposition}\label{subgroup}
\mbox{}
\begin{enumerate}\sep
\item If $R$ is a ring or an $F$-algebra then $N(R)$ is closed under conjugation and inversion. If $R$ is commutative then $N(R)$ is a subgroup of $Q(R)$.
\item If $R$ is an $F$-algebra, then $\pi_F(R)$ is closed under conjugation and inversion. If $R$ is commutative then $\pi_F(R)$ is a subgroup of $Q(R)$.
\item If $R$ is a ring, then $\pi(R)$ is closed under conjugation. If $R$ is commutative then $\pi(R)$ is a submonoid of $Q(R)$.
\end{enumerate}
\end{proposition}

\proof Let $a \in R$, $r \in Q(R)$, and let $p$ be a polynomial. Then by Lemma~\ref{aut} $r\circ p(a)\circ r^{(-1)}=p(r\circ a\circ r^{(-1)})$, so $r\circ a\circ r^{(-1)}$ is annihilated by the same polynomials as $a$. This shows that $N(R)$ and $\pi(R)$ (resp. $\pi_F(R)$) are closed under conjugation.

The group inverse (quasi-inverse) of a nilpotent element is a polynomial in this element, so it is again nilpotent. Thus $N(R)$ is closed under inversion. If $R$ is commutative then $N(R)$ is a subring of $R$, so it is closed under $\circ$ as well.

Let $R$ be an $F$-algebra. If $R$ is commutative then $A_F(R)$ is a subalgebra of $R$. Thus $A_F(R)$ is closed under $\circ$ and by Lemma~\ref{inclusion} so is $\pi_F(R)$.
If $a \in \pi_F(R)$ then the quasi-inverse of $a$ is a polynomial in $a$, so it is algebraic and hence contained in $A_F(R) \cap Q(R)=\pi_F(R)$.

If $R$ is a ring then by Proposition~\ref{integral} $\pi(R)^{(-1)}=I(R)\cap Q(R)$. If $R$ is commutative then $I(R)$ is a subring of $R$ and hence closed under $\circ$. So $\pi(R)^{(-1)}$ and consequently $\pi(R)$ is closed under $\circ$. \endproof

For a ring $R$ the set $\pi(R)$ need not be closed under inversion. For example the quasi-inverse of $1+\frac{1}{2} \in \pi(\QQ)$ is $1+2$ and is not contained in $\pi(\QQ)$. In fact, we know that $\pi(R)^{(-1)}=I(R)\cap Q(R)$.

\begin{example}\label{notcl}
Let $F$ be an algebraically closed field and $E=F(x)$ the field of rational functions over $F$. By Example~\ref{mat}, $\pi_F(M_2(E))$ consists of matrices with eigenvalues in $F \backslash\{1\}$. Take matrices
$$A=\left[\begin{array}{cc}0&x\\0&0 \end{array}\right] \qquad\textrm{and}\qquad B=\left[\begin{array}{cc}0&0\\1&0 \end{array}\right],$$
which both lie in $\pi_F(M_2(E))$, since they are nilpotent. Then
$$A \circ B=\left[\begin{array}{cc}-x&x\\1&0 \end{array}\right]$$
does not have eigenvalues in $F$, since its trace is $-x \notin F$. So $\pi_F(M_2(E))$ is not closed under $\circ$.
\end{example}

For a subset $S$ of $Q(R)$ let $\gen{S}$ denote the normal subgroup of $Q(R)$ generated by $S$. By Proposition~\ref{subgroup} we have:
\begin{eqnarray*}
\gen{N(R)} &=& \textup{ finite products of elements of } N(R),\\
\gen{\pi_F(R)} &=& \textup{ finite products of elements of } \pi_F(R),\\
\gen{\pi(R)} &=& \textup{ finite products of elements of } \pi(R) \cup \pi(R)^{(-1)},\\
\gen{\pi(R)\cap\pi(R)^{(-1)}} &=& \textup{ finite products of elements of } \pi(R)\cap I(R),
\end{eqnarray*}
where products means products in operation $\circ$.

\begin{example}
From Example~\ref{rac} it is easy to calculate that we have $\gen{\pi(\QQ)}=Q(\QQ)=\QQ \backslash\{1\}$ and $\gen{\pi(\QQ)\cap\pi(\QQ)^{(-1)}}=\{0,2\}$.
\end{example}

\begin{example}
Recall that a complex matrix $A$ is called \textit{unipotent} if $I-A$ is nilpotent, where $I$ denotes the identity matrix. In \cite{Wan} it was shown that a complex matrix is a finite product of unipotent matrices iff it has determinant $1$. This shows that
$$\gen{N(M_n(\CC))}=\set{A \in M_n(\CC)}{\det(I-A)=1},$$
which is a proper subgroup of
$$\pi_\CC(M_n(\CC))=Q(M_n(\CC))=\set{A \in M_n(\CC)}{\det(I-A) \neq 0}.$$
\end{example}

Next we investigate what can be said about addition. We will need the following proposition which may be of independent interest. Recall that an integral domain $K$ is called a \emph{factorization domain} (also an \emph{atomic domain}) if every nonzero nonunit of $K$ can be written as a finite product of irreducible elements.

\begin{proposition}\label{div}
Let $R$ be a unital ring and $K$ a commutative subring of $R$ with $1 \in K$ such that $R \backslash K \subseteq R^{-1}$. If $K$ is a factorization domain then one of the following holds:
\begin{enumerate}\sep
\item $R=K$,
\item $R$ is a local ring with maximal ideal $m \subseteq K$ and $K$ is a local ring with maximal ideal $m$,
\item $R$ is a division ring.
\end{enumerate}
\end{proposition}

\proof Suppose that $R \neq K$ and $R$ is not a division ring. Then there exist $r \in R\backslash K \subseteq R^{-1}$ and $0 \neq a \in K \backslash R^{-1}$. Since $K$ is a factorization domain, we may assume that the element $a$ is irreducible. We will prove that $K^{-1}=K\cap R^{-1}$. Let $x$ be arbitrary element of $K$ that is invertible in $R$ and set $y=x^{-1}a$. Then $y$ is not invertible in $R$, since $a$ is not. But $R\backslash K \subseteq R^{-1}$, so $y \in K$. Thus $a=xy$ is a factorization of $a$ in $K$. Since $a$ was irreducible and $y$ is not invertible, $x$ must be invertible in $K$, as needed. Now let $m$ be the set of all elements of $K$ that are not invertible in $K$. Since $R \backslash K \subseteq R^{-1}$, $m$ is also the set of all non-invertible elements of $R$. If $x \in m$ and $k \in K$ then $xk$ is not invertible in $K$, otherwise $x$ would be invertible due to the commutativity of $K$. So $mK \subseteq m$. If $x,y \in m$ then by the above $x$ and $y$ are not invertible in $R$. By the choice of $r$ this implies that $xr$ and $yr$ are not invertible in $R$, so $xr,yr \in K$. Thus $(x-y)r \in K$. But $x-y \in K$ and $r \notin K$, hence $x-y$ cannot be invertible in $K$, so $x-y \in m$. This proves that $m$ in an ideal in $K$, so $K$ is local with maximal ideal $m$. Now let $x \in m$ and $s \in R$, so by the above $x$ is not invertible in $R$. If $s \in K$ then $sx,xs \in m$ by what we have just proved. If $s \notin K$ then $s$ is invertible in $R$. So $sx$ and $xs$ are not invertible in $R$, hence $sx,xs \in m$. This shows that $m$ is also an ideal of $R$ and $R$ is local with maximal ideal $m$. \endproof

\begin{remark}
There exist examples where case (ii) of Proposition~\ref{div} occurs in a nontrivial way. Take for example $R=E[[x]]$ and $K=F+E[[x]]x \subseteq R$ where $F \varsubsetneq E$ are fields. Every nonzero nonunit in $K$ is contained in $E[[x]]x$ and factors as $x^ng(x)$ for some nonnegative integer $n$ and some $g(x)$ of the form $\alpha_1 x+\alpha_2x^2+\alpha_3x^3+\ldots$ with $\alpha_1 \neq 0$.
\end{remark}

\begin{theorem}\label{add}
Let $R$ be a ring. For any subgroup $S$ of $Q(R)$ the following are equivalent:
\begin{enumerate}\sep
\item $S$ is closed under addition,
\item $S$ is closed under multiplication,
\item $S$ is a subring of $R$.
\end{enumerate}
\end{theorem}

\proof We can verify by a short calculation that for any $x,y \in Q(R)$ we have
$$xy=x\circ (x^{(-1)}+y^{(-1)})\circ y \qquad\textup{and}\qquad x+y=x\circ (x^{(-1)}y^{(-1)})\circ y.$$
This shows that (i) and (ii) are equivalent. For $x \in Q(R)$ we also have
$$-x=(2x^{(-1)})\circ x,$$
so (i) implies that $S$ is closed under negation as well, which implies (iii). Clearly (iii) implies (i). \endproof

As a corollary to Theorem~\ref{add} we have the following.

\begin{corollary}\label{cor1}
Let $F$ be a field of characteristic $0$ and $R$ a commutative $F$-algebra. If $\pi_F(R)$ is closed under addition then $\pi_F(R)=N(R)$.
\end{corollary}

\proof Since $R$ is commutative, $\pi_F(R)$ is a subgroup of $Q(R)$ by Proposition~\ref{subgroup}. If $\pi_F(R)$ is closed under addition then it is a subring of $R$ by Theorem~\ref{add}. Let $a \in R$ be $\pi$-algebraic with polynomial $p$ and let $\lambda$ be a nonzero scalar. Since $F$ is of characteristic $0$ there exists a positive integer $n$ such that $n\lambda^{-1}$ is not a zero of $p$. Hence $n^{-1}\lambda a$ is $\pi$-algebraic with polynomial $p(n\lambda^{-1})^{-1}p(n\lambda^{-1}x)$. Since $\pi_F(R)$ is closed under addition and $\lambda a$ is a multiple of $n^{-1}\lambda a$, $\lambda a$ is $\pi$-algebraic as well. So $\pi_F(R)$ is in fact a subalgebra of $R$. Thus $\pi_F(R)$ is nil by Proposition~\ref{alg} and $\pi_F(R)=N(R)$ follows. \endproof

The conclusion of Corollary~\ref{cor1} also holds for rings.

\begin{proposition}\label{pro1}
Let $R$ be a commutative ring. If $\pi(R)$ is closed under addition then $\pi(R)=N(R)$.
\end{proposition}

\proof Suppose $\pi(R)$ is closed under addition. First we show that $\pi(R)$ is closed also under negation. If $a$ is $\pi$-algebraic, then $\NN a \subseteq \pi(R)$ since $\pi(R)$ is closed under addition. By Proposition~\ref{na} and Remark~\ref{rem-na} there exists a nonzero integer $n$ such that $na$ is nilpotent. Thus $-|n|a$ is nilpotent and hence $\pi$-algebraic. So $-a=-|n|a+(|n|-1)a$ is $\pi$-algebraic as well, since $(|n|-1)a$ is a nonnegative multiple of $a$. The commutativity of $R$ implies that $\pi(R)$ is closed under $\circ$. Since $xy=x+y-x\circ y$, $\pi(R)$ is closed under multiplication as well. So $\pi(R)$ is a $\pi$-algebraic subring of $R$, hence it is nil by Theorem~\ref{ring}. \endproof

We are now left with the case of algebras over fields of prime characteristic. We were not able to obtain an analogue of Corollary~\ref{cor1} for arbitrary fields of prime characteristic, but only for algebraic extensions of prime fields.

\begin{corollary}\label{cor2}
Let $p$ be a prime number, $F$ an algebraic field extension of the prime field $\ZZ/p\ZZ$, and $R$ a commutative $F$-algebra. If $\pi_F(R)$ is closed under addition then $\pi_F(R)=N(R)$.
\end{corollary}

\proof Since $F$ is algebraic over $\ZZ/p\ZZ$, we have $A_F(R)=A_{\ZZ/p\ZZ}(R)$, so $\pi_F(R)=\pi_{\ZZ/p\ZZ}(R)$ by Lemma~\ref{inclusion}. Now let $a \in R$ be $\pi$-algebraic over $\ZZ/p\ZZ$ with polynomial $\widehat{f}$ and let $f$ be a polynomial with integer coefficients that represents $\widehat{f}$. Since $\widehat{f}(1)=1$, there exists an integer $k$ such that $f(1)=kp+1$. If we set $F(x)=f(x)-kpx$, then $F(0)=0$, $F(1)=1$ and $F(a)=0$, since $pa=0$. So $a$ is $\pi$ algebraic over $\ZZ$. Hence $\pi_{\ZZ/p\ZZ}(R) \subseteq \pi(R)$ and clearly $\pi(R) \subseteq \pi_{\ZZ/p\ZZ}(R)$. This implies $\pi_F(R)=\pi(R)$ and so $\pi_F(R)=N(R)$ by Proposition~\ref{pro1}. \endproof

This was one extremal situation, when every $\pi$-algebraic element is in fact nilpotent. The other extremal situation would be when there are no nilpotent elements, but many $\pi$-algebraic ones. As we have mentioned before, in an algebraic division algebra there are no nonzero nilpotent elements although all elements except the unit are $\pi$-algebraic. Next we investigate when something similar happens in general algebras. The question is whether $\pi_F(R) \cup \{1\}$ will form a division subring of a unital $F$-algebra $R$. When $R$ is just a ring, we can ask a similar question, however it seems more natural to consider the set $\gen{\pi(R)} \cup (\ZZ\cdot 1)$ in this case, since the elements in $(\ZZ\cdot 1) \backslash \{1\}$ need not be automatically contained in $\gen{\pi(R)}$. In certain situations though, they are.

\begin{theorem}\label{addunitzero}
Let $R$ be a unital ring of characteristic $0$. For any subgroup $S$ of $Q(R)$ with $\{0,2\} \varsubsetneq S$ the following are equivalent:
\begin{enumerate}\sep
\item\label{auz1} $S \cup \ZZ$ is closed under addition,
\item\label{auz2} $S \cup \{1\}$ is a division subring of $R$.
\end{enumerate}
\end{theorem}

\proof Clearly \eqref{auz2} implies \eqref{auz1}, since in this case $S \cup \{1\}=S \cup \ZZ$. So assume \eqref{auz1} holds. First we show that $S \cup \ZZ$ is a subring. If $x \in S \cup \ZZ$ then $2 \circ x=2-x \in S \cup \ZZ$, since $2 \in S \cap \ZZ$. So if $x \in S \cup \ZZ$ then $-x=2-(2+x) \in S \cup \ZZ$ by \eqref{auz1}. Thus $S \cup \ZZ$ is closed under negation. $S$ and $\ZZ$ are both closed under $\circ$. If $x \in S$ and $n \in \ZZ$ then $x \circ n=n \circ x=n+x-nx  \in S \cup \ZZ$, since $nx$ is a multiple of $x$ or $-x$ and $S \cup \ZZ$ is closed under addition. So $S \cup \ZZ$ is closed under $\circ$ and also under multiplication since $xy=x+y-x\circ y$. This shows that $S \cup \ZZ$ is a subring of $R$. Now every element in $S$ is quasi-regular with quasi-inverse in $S$, thus every element in $1-S$ is invertible in $S \cup \ZZ$. Since $S \cup \ZZ$ is a subring, we have $1-S \backslash \ZZ=S \backslash \ZZ$. So every element in $S \backslash \ZZ$ is invertible in $S \cup \ZZ$. By Proposition~\ref{div} either $S \subseteq \ZZ$ or $S \cup \ZZ$ is a division ring. Suppose $S \subseteq \ZZ$. Then the quasi-inverse of every element in $S \subseteq \ZZ$ lies again in $S \subseteq \ZZ$, so $S\subseteq Q(\ZZ)=\{0,2\}$, which contradicts our assumption. Therefore $S \cup \ZZ$ is a division ring. It remains to prove that $\ZZ \backslash \{1\} \subseteq S$. Let $n \in \ZZ \backslash \{1\}$. If $n=0$ or $n=2$ then $n \in S$ by assumption. So suppose $n \neq 0,2$. Since $S \cup \ZZ$ is a division ring, $1-n$ is invertible in $S \cup \ZZ$. Since $1-n \neq \pm 1$, the fact that the characteristic of $R$ is $0$ implies $(1-n)^{-1} \notin \ZZ$, i.e. $1-(1-n)^{-1} \in S$.  Consequently $n=(1-(1-n)^{-1})^{(-1)} \in S$, since $S$ is a subgroup of $Q(R)$.
\endproof

\begin{theorem}\label{addunitprime}
Let $R$ be a unital ring of prime characteristic $p$. For any subgroup $S$ of $Q(R)$ the following are equivalent:
\begin{enumerate}\sep
\item $S \cup \ZZ/p\ZZ$ is closed under addition,
\item $S \cup \ZZ/p\ZZ$ is a division subring of $R$.
\end{enumerate}
\end{theorem}

\proof In this case $S \cup \ZZ/p\ZZ$ is automatically closed under negation, since $-x=(p-1)x$ is a multiple of $x$. The proof is now the same as that of Theorem~\ref{addunitzero} except for the case $S \subseteq \ZZ/p\ZZ$, but in this case $S \cup \ZZ/p\ZZ=\ZZ/p\ZZ$ is automatically a division ring. \endproof

\begin{corollary}
Let $F$ be a field and $R$ a unital commutative $F$-algebra. If $\pi_F(R) \cup \{1\}$ is closed under addition then it is a subfield of $R$.
\end{corollary}

\proof This follows directly from Proposition~\ref{subgroup} and Theorems~\ref{addunitzero} and \ref{addunitprime}, since $(\ZZ \cdot 1) \backslash \{1\} \subseteq \pi_F(R)$ by Lemma~\ref{2minus}. \endproof

\begin{corollary}
Let $R$ be a unital commutative ring of prime or $0$ characteristic with $\pi(R) \neq \{0,2\}$. If $\pi(R)\circ \pi(R)^{(-1)} \cup (\ZZ \cdot 1)$ is closed under addition then it is a subfield of $R$.
\end{corollary}

\proof The commutativity of $R$ implies $\langle \pi(R)\rangle=\pi(R)\circ \pi(R)^{(-1)}$. Since $2 \in \pi(R)$, the result follows from Theorems~\ref{addunitzero} and \ref{addunitprime}. \endproof

\end{document}